\documentclass[12pt]{amsart}
\usepackage[utf8]{inputenc}
\usepackage[english]{babel}
\usepackage{tikz-cd}
\usepackage{verbatim}

\title{Manin's \MakeLowercase{b}-constant in families}
\author{Akash Kumar Sengupta}
\date{}

\usepackage[
  hmarginratio={1:1},     % equal left and right margins
  vmarginratio={1:1},     % equal top and bottom margins
  textwidth=410pt,        % new text width
  heightrounded,          % always useful
]{geometry}
\usepackage{lipsum}%opening
\usepackage{amsmath,amsthm,amssymb,mathrsfs}

\newtheoremstyle{common}
    {6pt plus 5pt minus 2pt}% above space (default)
    {6pt plus 5pt minus 2pt}% below space
    {\normalfont}% body font
    {0em}% indent
    {\bfseries}% head font
    {}% punct after head
    {.5em}% space
    {}% custom
\theoremstyle{common}

\newtheorem{thm}{Theorem}[section]
\newtheorem{lem}[thm]{Lemma}
\newtheorem{prop}[thm]{Proposition}
\newtheorem{rem}[thm]{Remark}
\newtheorem{defn}[thm]{Definition}

\begin{document}

\maketitle
\begin{abstract}
    We show that the $b$-constant (appearing in Manin's conjecture) is constant on very general fibers of a family of algebraic varieties. If the fibers of the family are uniruled, then we show that the $b$-constant is constant on general fibers.
\end{abstract}
\section{Introduction}
Let $X$ be a smooth projective variety over a field of $k$ of characteristic $0$ and $L$ a big Cartier $\mathbb{Q}$-divisor on $X$. Let $\Lambda_{\mathrm{eff}}(X) \subset \mathrm{NS}(X)_{\mathbb{R}}$ be the cone of pseudo-effective divisors. The Fujita invariant or the $a$-constant is defined as  
\[a(X,L)= \mathrm{min}\{t\in \mathbb{R}| [K_X]+t[L] \in \Lambda_{\mathrm{eff}}(X)\}\]
The invariant $\kappa\epsilon(X,L)=-a(X,L)$ was introduced and studied by Fujita under the name Kodaira energy in [Fuj87], [Fuj92]. The $a$-constant was introduced in the context of Manin's conjecture in [FMT89].\\

The $b$-constant is defined as follows (cf. [FMT89], [BM90])
\[ b(X,L) = \mathrm{codim \hspace*{0.1cm} of \hspace*{0.1cm}minimal\hspace*{0.1cm} supported\hspace*{0.1cm} face\hspace*{0.1cm} of }\hspace*{0.1cm} \Lambda_{\mathrm{eff}}(X)\hspace*{0.1cm}\] \[\mathrm{containing \hspace*{0.1cm} the\hspace*{0.1cm} class\hspace*{0.1cm} of\hspace*{0.1cm} }K_X+a(X,L)L\]
For a singular variety $X$, the $a$,$b$-constants of $L$ are defined to be the $a$,$b$-constants of $\pi^*L$ on a resolution $\pi:\widetilde{X}\longrightarrow X$. \\

Let $f: X\longrightarrow T$ be a family of projective varieties and $L$ a $f$-big and $f$-nef Cartier $\mathbb{Q}$-divisor. By semi-continuity the $a$-constant of the fibers $a(X_t,L|_{X_t})$  is constant on very general fiber. It follows from invariance of log plurigenera that if the fibers are uniruled then the $a$-constant is constant on general fibers.\\

In this paper we investigate the behaviour of the $b$-constant in families and answer the questions posed in [LT17a]. We prove the following:

\begin{thm}
Let $f:X\longrightarrow T$ be a projective morphism of irreducible varieties over an algebraically closed field $k$ of characteristic $0$, such that the generic fiber is geometrically integral. Let $L$ be a $f$-big Cartier $\mathbb{Q}$-divisor. Then there exists a countable union of proper closed subvarieties $Z= \cup_i Z_i\subsetneq T$, such that
\[b(X_{\overline{t}},L|_{X_{\overline{t}}})= b(X_{\overline{\eta}},L|_{X_{\overline{\eta}}})\]
for all $t\in T\setminus Z$, where $\eta \in T$ is the generic point. In particular, the $b$-constant is constant on very general fibers.\\
\end{thm}

If the fibers of the family are uniruled, then we have the following:

\begin{thm}
Let $f:X\longrightarrow T$ be a projective morphism of irreducible varieties over an algebraically closed field $k$ of characteristic $0$, such that the generic fiber is geometrically integral. Let $L$ be a $f$-big and $f$-nef Cartier $\mathbb{Q}$-divisor. Suppose the fibers $X_t$ are uniruled for all $t \in T$. Then there exists a proper closed subscheme $W\subsetneq  T$ such that \[b(X_{\overline{t}},L|_{X_{\overline{t}}})= b(X_{\overline{\eta}},L|_{X_{\overline{\eta}}})\]
for $t\in T\setminus W$ and $\eta \in T$ is the generic point. In particular, the $b$-constant is constant on general fibers in a family of uniruled varieties.\\
\end{thm}

One can not replace the very general condition in Theorem $1.1$ by just general. For example, in a family of K$3$-surfaces the $b$-constant of a fiber is same as the Picard rank and there exist families where the Picard rank jumps on infinitely many subvarities. Invariance of the $b$-constant in general fiber of a family of uniruled varieties was proved in [LT17a] under the assumption $\kappa(K_{X_t}+a(X_t,L|_{X_t})L|_{X_t})=0$. Theorem $1.2$ generalizes their result to get rid of this condition on fibers.\\

One of the motivations for studying the behaviour of $a$ and $b$-constants is Manin's conjecture about asymptotic growth of rational points on Fano varieties proposed in [FMT89], [BM90]. The following version was suggested by Peyre in [Pey03] and later stated in [Rud14], [BL15].\\

{\bf Manin's conjecture:} Let $X$ be a Fano variety defined over a number field $F$ and $\mathcal{L}=(L,||.||)$ a big and nef adelically metrized line bundle on $X$ with associated height function $H_{\mathcal{L}}$. Then there exists a thin set $Z \subset X(F)$ such that one has 
\[\#\{x\in X(F)-Z| H_{\mathcal{L}}(x)\leq B\} \sim c(F,X(F)-Z,\mathcal{L})B^{a(X,L)}\mathrm{log}B^{b(X,L)-1}\]
as $B\rightarrow \infty$. \\

For the geometric consistency of Manin's conjecture, a necessary condition is that $a$,$b$-constants acheive a maximum as we vary over subvarieties of $X$. The behaviour of the $a$,$b$-constants in families was used in [LT17a] to show this necessary condition. The $a$ and $b$-constants also play a role in determining and counting the dominant components of the space $\mathrm{Mor}(\mathbb{P}^1,X)$ of morphisms from $\mathbb{P}^1$ to a smooth Fano variety $X$ (see [LT17b] for details). \\

The ideas in proving our results are as follows. To prove Theorem 1.1, we analyze the behaviour of the $b$-constant under specialization and combine this with the constancy of the Picard rank and the $a$-constant in very general fibers to obtain the desired conclusion. The key step for Theorem 1.2 is to  prove constancy on closed points when $k=\mathbb{C}$. We run a $K_X+aL$-MMP over the base $T$, to obtain a relative minimal model $X\dashrightarrow X'$ where $a=a(X_t,L|_{X_t})$. We pass to a relative canonical model
$ \phi: X\dashrightarrow Z$ over $T$ and base change to $t\in T$, to obtain $\phi_t:X_t \dashrightarrow Z_t$ as the canonical model for $(X_t,aL_{X_t})$. Using a version of the global invariant cycles theorem (see Lemma 2.11), we observe that $b(X_t,L_t)$ is same as the rank of the mondromy invariant subspace of $N^1(Y'_z)_{\mathbb{R}}$, where $Y'_z$ is a general fiber of $X'_t\longrightarrow Z_t$. Then using topological local triviality of algebraic morphisms we conclude that the monodromy invariant subspace has constant rank.\\

The outline of the paper is as follows. In Section $2$ we discuss the preliminaries. In Section $3$ and $4$ we prove Theorems $1.1$ and $1.2$ respectively.\\

{\bf Acknowledgements:} I am very grateful to my advisor Professor J\'{a}nos Koll\'{a}r for fruitful discussions and his constant support and encouragement. I am thankful to Sebasti\'{a}n Olano and Ziquan Zhuang for interesting conversations.\\

\section{Preliminaries}
In this paper we always work in characteristic $0$.

\subsection{N\'{e}ron-Severi group.} Let $X$ be a smooth proper variety over a field $k$. The N\'{e}ron-Severi group $\mathrm{NS}(X)$ is defined as the quotient of the group of Weil divisors, $\mathrm{Wdiv}(X)$, modulo algebraic equivalence. We denote $N^1(X)= \mathrm{Cdiv}(X)\slash \equiv$, the quotient of Cartier divisors by numerical equivalence. We denote $\mathrm{NS}(X)_{\mathbb{R}}= \mathrm{NS}(X) \otimes \mathbb{R}$ and similarly $N^1(X)_{\mathbb{R}}$.
By [N\'{e}r52], $\mathrm{NS}(X)_{\mathbb{R}}$ is a finite-dimensional vector space and its rank $\rho(X)$ is called the Picard rank. If $X$ is a smooth projective variety, then $ \mathrm{NS}(X)_{\mathbb{R}} \cong N^1(X)_{\mathbb{R}}$. 

\begin{rem}
Let $X$ be a smooth variety over an algebraically closed field $k$. If $k \subset k'$ is an extension of algebraically closed fields, then the natural homomorphism $\mathrm{NS}(X) \longrightarrow NS(X_{k'})$ is an isomorphism. So the Picard rank is unchanged under base extension of algebraically closed fields. 
\end{rem}

Let $X\longrightarrow T$ be a smooth proper morphism of irreducible varieties. Suppose $s$,$t\in T$ are such that $s$ is a specialization of $t$, i.e. $s$ is in the closure of $\{t\}$. Let $X_{\overline{t}}$ denote the base change to the algebraic closure of the residue field $k(t)$.

\begin{prop} (see [MP12, Prop. 3.6.]) In the situation as above, it is possible to choose a specialization homomorphism  
\[\mathrm{sp}_{\overline{t},\overline{s}}: \mathrm{NS}(X_{\overline{t}})\longrightarrow \mathrm{NS}(X_{\overline{s}})\]
such that 
\begin{itemize}
    \item [(a)] $\mathrm{sp}_{\overline{t},\overline{s}}$ is injective. In particular $\rho(X_{\overline{s}})\geq \rho(X_{\overline{t}})$.
    \item [(b)] If $\mathrm{sp}_{\overline{t},\overline{s}}$ maps a class $[L]$ to an ample class, then $L$ is ample.

\end{itemize}
 
\end{prop}
If $\rho(X_{\overline{s}})= \rho(X_{\overline{t}})$, then the homomorphism $\mathrm{NS}(X_{\overline{t}})_{\mathbb{R}} \longrightarrow \mathrm{NS}(X_{\overline{s}})_{\mathbb{R}}$ is an isomorphism. \\

Let $X \longrightarrow T$ be a smooth projective morphism of irreducible varieties over $\mathbb{C}$. In Section 12 of [KM92], the local system $\mathcal{GN}^1(X/T)$ was introduced. This is a sheaf in the analytic topology defined as:
\[\mathcal{GN}^1(X/T)(U)=\{\text{sections of } \mathcal{N}^1(X/T) \text{ over } U \text{ with open support}\}\]
for analytic open $U\subset T$, and the functor $\mathcal{N}^1(X/T)$ is defined as $N^1(X\times_T T')$ for any $T'\longrightarrow T$. It was shown in [KM92, 12.2] that $\mathcal{GN}^1(X/T)$ is a local system with finite monodromy and $\mathcal{GN}^1(X/T)|_t=N^1(X_t)$ for very general $t\in T$. We can base change to a finite \'{e}tale cover of $T'\longrightarrow T$ so that $\mathcal{GN}^1(X'/T')$ has trivial monodromy. Then we have a natural identification between the fibers of $\mathcal{GN}^1(X'/T')$ and $N^1(X'/T')$. Therefore, for $t'\in T'$ very general, the natural map $N^1(X'/T') \longrightarrow N^1(X'_{t'})$ is an isomorphism. One can prove the same results over any algebraically closed field of characteristic $0$, by using the Lefschetz principle.\\

\subsection{Geometric invariants.} The pseudo-effective cone $\Lambda_{\mathrm{eff}}(X)$ is the closure of the cone of effective divisor classes in $\mathrm{NS}(X)_{\mathbb{R}}$. The interior of $\Lambda_{\mathrm{eff}}(X)$ is the cone of big divisors $\mathrm{Big}^1(X)_{\mathbb{R}}$. \\

\begin{defn} Let $L$ be a big Cartier $\mathbb{Q}$ divisor on $X$. The $a$-constant is 
\[a(X,L)= \mathrm{min}\{t\in \mathbb{R}| K_X+tL \in \Lambda_{\mathrm{eff}}(X)\}\]

\end{defn}
 For a singular projective variety we define $a(X,L):= a(\widetilde{X},\pi^*L)$ where $\pi: \widetilde{X}\longrightarrow X$ is a resolution of $X$. It is invariant under pull-back by a birational morphism of smooth varieties and hence independent of the choice of the resolution. By [BDPP13] we know that  $a(X,L)>0$ if and only if $X$ is uniruled. We note that, by flat base change, the $a$-constant is independent of base change to another field. \\

It was shown in [BCHM10] that, if $X$ is an uniruled with klt singularities and $L$ is ample, then $a(X,L)$ is a rational number. If $L$ is big and not ample, then $a(X,L)$ can be irrational (see [HTT15, Example 6]). For a smooth projective variety $X$, the function $a(X,\_):\mathrm{Big}^1(X)_{\mathbb{R}} \longrightarrow \mathbb{R}$ is a continuous function (see [LTT16, Lemma 3.2]).\\

\begin{defn}
A morphism $f: X \longrightarrow T$ between irreducible varieties is called a family of varieties if the generic fiber is geometrically integral. A family of projective varieties is a projective morphism which is a family of varieties.\\
\end{defn}

We recall the following result about the $a$-constant in families:
\begin{thm}([LT17a], [HMX13]).
Let $f: X \longrightarrow T$ be a family of uniruled projective varieties over an algebraically closed field. Let $L$ a $f$-big and $f$-nef $\mathbb{Q}$-Cartier divisor on $X$. Then there exists a nonempty subset $U\subset T$ such that $a(X_t,L|_{X_t})$ is constant for $t\in U$ and the Iitaka dimension $\kappa(K_{X_t}+a(X_t,L|_{X_t})L|_{X_t})$ is constant for $t\in U$.\\

\end{thm}

\begin{defn}
Let $X$ be a smooth projective variety over $k$ and $L$ a big Cartier $\mathbb{Q}$-divisor. The $b$-constant is defined as
\[ b(k,X,L) = \mathrm{codim \hspace*{0.1cm} of \hspace*{0.1cm}minimal\hspace*{0.1cm} supported\hspace*{0.1cm} face\hspace*{0.1cm} of }\hspace*{0.1cm} \Lambda_{\mathrm{eff}}(X)\hspace*{0.1cm}\] \[\mathrm{containing \hspace*{0.1cm} the\hspace*{0.1cm} class\hspace*{0.1cm} of\hspace*{0.1cm} }K_X+a(X,L)L\]
\end{defn}

It is invariant under pullback by a birational morphism of smooth varieties ([HTT15]). For a singular variety $X$ we define $b(k,X,L):=b(k,\widetilde{X},\pi^*L)$, by pulling back to a resolution. By Remark 2.1, if we have an extension $k \subset k'$ of algebraically closed fields, the pull back map $\mathrm{NS}(X) \longrightarrow NS(X_{k'})$ is an isomorphism and the pseudo-effective cones are isomorphic by flat base change. Also, $K_X+a(X,L)L$ maps to $K_{X_{k'}}+a(X_{k'},L_{k'})L_{k'}$ under this isomorphism. Therefore the $b$-constant is unchanged, i.e. $b(k',X_{k'},L_{k'})=b(k,X,L)$. From now on, when our base field is algebraically closed we write $b(X,L)$ instead of $b(k,X,L)$.\\

\subsection{Minimal and Canonical models.} Let $(X,\Delta)$ be a klt pair, with $\Delta$ a $\mathbb{R}$-divisor and $K_X+\Delta$ is $\mathbb{R}$-Cartier. Let $f: X \longrightarrow T$ be a projective morphism. A pair $(X',\Delta')$ sitting in a diagram
\[\begin{tikzcd}
X \arrow[rr, dashrightarrow, "\phi"] \arrow[dr, "f"]& & X'\arrow[dl,"f'"]\\
& T &\\
\end{tikzcd}\]
is called a $\mathbb{Q}$-factorial minimal model of $(X,\Delta)$ over $T$ if 
\begin{itemize}
    \item [(1)] $X'$ is $\mathbb{Q}$-factorial,
    \item[(2)] $f'$ is projective,
    \item[(3)] $\phi$ is a birational contraction,
    \item[(4)]$\Delta'=\phi_*\Delta$ and $(X',\Delta')$ is dlt,
    \item[(5)] $K_{X'}+\Delta'$ is $f'$-nef,
    \item[(6)] $a(E,X,\Delta) < a(E,X',\Delta')$ for all $\phi$-exceptional divisors $E \subset X$. Equivalently, if for a common resolution $p: W \longrightarrow X$ and $q: W \longrightarrow X'$, we may write 
    \[p^*(K_X+\Delta) = q^* (K_{X'}+\Delta') +E\]
    where $E\geq 0$ is $q$-exceptional and the support of $E$ contains the strict transform of the $\phi$-exceptional divisors.

\end{itemize}

A canonical model over $T$ is defined to be a projective morphism $g: Z \longrightarrow T$ with a surjective morphism $\pi: X' \longrightarrow Z$ with connected geometric fibers from a minimal model such that $K_{X'}+\Delta'=\pi^*H$ for an $\mathbb{R}$-Cartier divisor $H$ on $Z$ which is ample over $T$.\\

Suppose $K_X+\Delta$ is $f$-pseudo-effective and $\Delta$ is $f$-big, then by [BCHM10], we may run a $K_X+\Delta$-MMP with scaling to obtain a $\mathbb{Q}$-factorial minimal model $(X',\Delta')$ over $T$. It follows that $(X',\Delta')$ is also klt. Then the basepoint freeness theorem implies that $(K_X+\Delta')$ is $f'$ semi-ample. Hence there exists a relative canonical model $g: Z\longrightarrow T$.  In particular, if $\Delta$ is a $\mathbb{Q}$-divisor, the $\mathcal{O}_T$-algebra
\[\mathfrak{R}(X',\Delta')= \oplus_m f'_*\mathcal{O}_{X'}(\lfloor m(K_{X'}+\Delta')\rfloor )\]
is finitely generated. Let $X' \longrightarrow Z \longrightarrow \mathcal{P}\mathrm{roj}(\mathfrak{R}(X',\Delta'))$ be the Stein factorization of the natural morphism. Then $Z$ is the relative canonical model over $T$.\\

The following result relates the relative MMP over a base to the MMP of the fibers ( see [dFH11, Theorem 4.1] and [KM92, 12.3] for related statements).

\begin{lem}
Let $f:X\longrightarrow T$ be a flat projective morphism of normal varieties. Suppose $X$ is $\mathbb{Q}$-factorial and $D$ be an effective $\mathbb{R}$-divisor such that $(X,D)$ is klt. Let $\psi: X\longrightarrow Z$ be the contraction of a $K_X+D$-negative extremal ray of $\overline{\mathrm{NE}}(X/T)$. Suppose for $t\in T$ very general, the restriction map $N^1(X/T)\longrightarrow N^1(X_t)$ is surjective and $X_t$ is $\mathbb{Q}$-factorial.\\

Let $t\in T$ be very general. If $\psi_t:X_t\longrightarrow Z_t$ is not an isomorphism, then it is a contraction of a $K_{X_t}+D_t$-negative extremal ray, and :
\begin{itemize}
    \item[(a)] If $\psi$ is of fiber type, so is $\psi_t$.
    \item[(b)] If $\psi$ is a divisorial contraction of a divisor $G$, then $\psi_t$ is a divisorial contraction of $G_t$ and $N^1(Z/T)\longrightarrow N^1(Z_t)$ is surjective.
    \item[(c)] If $\psi$ is a flipping contraction and $\psi^+:X^+\longrightarrow Z$ is the flip, then $\psi_t$ is a flipping contraction and $X_t^+$ is the flip of $\psi_t: X_t\longrightarrow Z_t$. Also, $N^1(X^+/T)\longrightarrow N^1(X_t^+)$ is surjective.\\
\end{itemize}
\end{lem}

\begin{proof}
 Since the natural restriction map $N^1(X/T)\longrightarrow N^1(X_t)$ is surjective for very general $t\in T$, any curve in $X_t$ that spans a $K_X+D$-negative extremal ray $R$ of $\overline{\mathrm{NE}}(X/T)$, also spans a $K_{X_t}+D_t$ negative extremal ray $R_t$ of $\overline{\mathrm{NE}}(X_t)$. For $t\in T$ general, the base change $Z_t$ is normal and the morphism $X_t\longrightarrow Z_t$ has connected fibers, hence ${\psi_t}_*\mathcal{O}_{X_t}=\mathcal{O}_{Z_t}$. Hence $\psi_t$ is the contraction of the ray $R_t$ for very general $t\in T$.
 
 If $\psi$ is of fiber type, then so is $\psi_t$ for general $t\in T$. Let us assume that $\psi$ is birational.
 
 Suppose $\psi$ is a divisorial contraction of a divisor $G$. Then all components of $G_t$ are contracted. By the injectivity of $N_1(X_t)\longrightarrow N_1(X/T)$, we see that $\psi_t$ is an extremal divisorial contraction of $G_t$ (and $G_t$ is irreducible). Since $X_t$ is $\mathbb{Q}$-factorial, we have the surjectivity of $N^1(Z/T)\longrightarrow N^1(Z_t)$.
 
 Suppose $\psi$ is a flipping contraction and $\phi: X \dashrightarrow X^+$ is the flip. For very general $t\in T$, $X_t \longrightarrow Z_t$ is a small birational contraction of the ray $R_t$. Also, $X^+_t\longrightarrow Z_t$ is also small birational and $K_{X^+_t}+(\phi_*D)_t$ is $\psi^+$-ample for $t\in T$ general. Therefore $\phi_t: X_t \dashrightarrow X^+_t$ is the flip. The surjectivity of $N^1(X^+/T)\longrightarrow N^1(X_t^+)$ follows from $\psi_t$ being an isomorphism in codimension one.

\end{proof}

The next Proposition allows us to compare minimal and canonical models over a base to those of a general fiber.\\

\begin{prop}
 Let $f: X \longrightarrow T$ be a smooth morphism. Suppose $\Delta$ is an $f$-big and $f$-nef $\mathbb{R}$-divisor such that $(X,\Delta)$ is simple normal crossings pair with klt singularities. Suppose the local system $\mathcal{GN}^1(X/T)$ has trivial monodromy. Let $\phi: X\dashrightarrow X'$ be the relative minimal model obtained by running a $K_X+\Delta$-MMP over $T$ and $\pi: X'\longrightarrow Z$ be the morphism to the canonical model over $T$. Then for a general $t\in T$,
\begin{itemize}
    \item[(1)] The base change $\phi_t: X_t\dashrightarrow X'_t$ is a $\mathbb{Q}$-factorial minimal model of $(X_t,\Delta_t)$,
    \item[(2)] Also, $\pi_t: X'_t\longrightarrow Z_t$ is the canonical model of $(X_t,\Delta_t)$.\\
\end{itemize}
\end{prop}

\begin{proof}
(1) Since $\mathcal{GN}^1(X/T)$ has trivial monodromy, the natural restriction morphism $N^1(X/T)\xrightarrow{\sim} N^1(X_t)$ is an isomorphism for $t\in T$ very general. Then Lemma 2.7 implies that, for very general $t\in T$, the base change $\phi_t:X_t \dashrightarrow X'_t$ is a composition of steps of the $K_{X_t}+\Delta_t$-MMP. In particular, $X'_t$ is $\mathbb{Q}$-factorial for a very general $t\in T$. The fibers $X'_t$ have terminal singularities, by [LTT16, Lemma 2.4]. Hence [KM92, 12.1.10] implies that there is a non-empty open $U\subset T$ such that $X'_t$ is $\mathbb{Q}$-factorial for $t\in U$. For a general $t\in T$, the conditions (2)-(6) in the definition of a minimal model follows easily. Therefore, $(X'_t,\Delta'_t)$ is a $\mathbb{Q}$-factorial minimal model of $(X_t,\Delta_t)$ for general $t\in T$.\\

(2) Let $g: Z\longrightarrow T$ be the relative canonical model.
 Now $Z$ is normal. Therefore, for a general $t\in T$, the base change $Z_t$ is normal and $X'_t\longrightarrow Z_t$ has geometrically connected fibers. Also, $K_{X'}+\Delta=g^*H$ where $H$ is a $\pi$-ample $\mathbb{R}$-Cartier divisor on $Z$. By adjunction, $K_{X'_t}+\Delta'_t$ is pull-back of an ample $\mathbb{R}$-Cartier divisor on $Z_t$. Hence, $X'_t\longrightarrow Z_t$ is the canonical model for general $t\in T$.\\
\end{proof}

Let $X$ be a smooth uniruled projective variety over an algberaically closed field and $L$ a big and nef $\mathbb{Q}$-divisor on $X$. The following result (contained in [LTT16]) gives a geometric interpretation of the $b$-constant.

\begin{prop}
Let $\phi: X \dashrightarrow X'$ be a $K_X+a(X,L)L$-minimal model. Then
\begin{itemize}
    \item [(1)]$b(X,L)=b(X',\phi_*L)$.
    \item[(2)]If $\kappa(K_X+a(X,L)L) =0$ then $b(X,L)= \mathrm{rk}N^1(X')_{\mathbb{R}}$.
    \item[(3)] If $\kappa(K_X+a(X,L)L) >0$ and $\pi: X'\longrightarrow Z$ is the morphism to the canonical model and $Y'$ is a general fiber of $\pi$. Then
    \[b(X,L)=\mathrm{rk}N^1(X')_{\mathbb{R}}-\mathrm{rk}N_{\pi}^1(X')_{\mathbb{R}}= \mathrm{rk}(\mathrm{im}(N^1(X')_{\mathbb{R}}\longrightarrow N^1(Y')_{\mathbb{R}}))\]
    where $N_{\pi}^1(X')_{\mathbb{R}}$ is the span of the $\pi$-vertical divisors and $N^1(X')_{\mathbb{R}}\longrightarrow N^1(Y')_{\mathbb{R}}$ is the restriction map.
\end{itemize} 

\end{prop}
\begin{proof}
Part (1) is the statement of Lemma 3.5 in [LTT16]. Part (2) follows from part (1). By abundance, $K_X+a(X,L)\phi_*L$ is semi-ample. Then  $\kappa(K_X+a(X,L)L) =0$ implies that $K_X+a(X,L)\phi_*L\equiv 0$. Hence, $b(X,L)=b(X',\phi_*L)=\mathrm{rk}N^1(X')_{\mathbb{R}}$. Part (3) follows from the proof of Theorem 4.5 in [LTT16].

\end{proof}

\begin{prop}(see [LT17a, Prop. 4.4])
Let $f:X \longrightarrow T$ be a family of projective varieties. Suppose $L$ is a $f$-big and $f$-nef Cartier divisor on $X$. Assume that for a general member $\kappa(K_{X_t}+a(X_t,L_t)L_t)=0$. Then $b(X_t,L_t)$ is constant on for general $t\in T$.\\ 
\end{prop}

\subsection{Global invariant cycles.} Let $\pi: X\longrightarrow Z$ be a morphism of complex algebraic varieties. Then, by Verdier's generalization of Ehresmann's theorem [Ver76, Corrolaire 5.1], there exists a Zariski open $U\subset Z$ such that $\pi^{-1}(U)\longrightarrow U$ is a topologically locally trivial fibration (in the analytic topology), i.e. every point $z\in U$ has a neighbourhood $N \subset U$ in the analytic topology, such that there is a fiber preserving homeomorphism

\[\begin{tikzcd}
\pi^{-1}(N) \arrow[rr, "\sim"] \arrow[dr]& & {N\times F} \arrow[dl] \\
& N &
\end{tikzcd}\]
where $F=\pi^{-1}(z)$. Consequently we have a monodromy action of $\pi_1(U,z)$ on the cohomology of the fiber $H^i(X_z,\mathbb{R})$.\\

The following result is an adaptation of Deligne's global invariant cycles theorem [Del71] to the case of singular varieties, which helps us to compute the $b$-constant.\\

\begin{lem}
 Let $\pi: X\longrightarrow Z$ be a morphism of normal projective varieties over $\mathbb{C}$ where $X$ is $\mathbb{Q}$-factorial. Let $\mu: \widetilde{X}\longrightarrow X$ be a resolution of singularities. Let $U\subset Z$ is a Zariski open subset such that $\pi \circ \mu$ is smooth over $U$ and $(\pi \circ \mu)^{-1}(U)\longrightarrow U$ and $\pi^{-1}(U)\longrightarrow U$ are topologically locally trivial fibrations (in the analytic topology). Suppose for general $z \in U$, the fiber  $X_z:=\pi^{-1}(z)$ is rationally connected with rational singularities. Then
 \[\mathrm{im}(N^1(X)_{\mathbb{R}} \longrightarrow N^1(X_z)_{\mathbb{R}}) \simeq H^2(X_z,\mathbb{R})^{\pi_1(U,z)}\]
 for general $z\in U$, where  $H^2(X_z,\mathbb{R})^{\pi_1(U,z)}$ is the monodromy invariant subspace. 
\end{lem}

Note that by generic smoothness and the discussion above, given any resolution of singularities $\mu: \widetilde{X}\longrightarrow X$, we may choose a Zariski open $U\subset Z$ such that $\pi\circ \mu$ is smooth over $U$ and $(\pi \circ \mu)^{-1}(U)\longrightarrow U$ and $\pi^{-1}(U)\longrightarrow U$ are topologically locally trivial fibrations.\\

\begin{proof}

Let $\widetilde{X}_z$ be the fiber of $\pi\circ \mu$ over $z$. For $z\in U$ general, $\mu_z: \widetilde{X}_z\longrightarrow X_z$ is a resolution of singularities. Since $X_z$ is rationally connected, $\mathbb{Q}$-linear equivalence and numerical equivalence of $\mathbb{Q}$-Cartier divisors coincide, i.e. $\mathrm{Pic}(X_z)_{\mathbb{Q}} \simeq N^1(X_z)_{\mathbb{Q}}$. We know $h^1(\widetilde{X}_z,\mathcal{O}_{\widetilde{X}_z})=h^2(\widetilde{X}_z,\mathcal{O}_{\widetilde{X}_z})=0$ since $\widetilde{X}_z$ is smooth rationally connected.  We also have $h^1(X_z,\mathcal{O}_{X_z})=h^2(X_z,\mathcal{O}_{X_z})=0$, because $X_z$ has rational singularities. Therefore $H^2(\widetilde{X}_z,\mathbb{Q})\simeq N^1(\widetilde{X}_z)_{\mathbb{Q}}$ and  $H^2(X_z,\mathbb{Q})\simeq N^1(X_z)_{\mathbb{Q}}$.\\

Consider the natural restriction map on cohomology groups $H^2(\widetilde{X},\mathbb{Q})\longrightarrow H^2(\widetilde{X}_z,\mathbb{Q})$. By Deligne's global invariant cycles theorem ([Del71] or [Voi03, 4.3.3]) we know that for $z\in U$,
\[\mathrm{im}(H^2(\widetilde{X},\mathbb{Q})\longrightarrow (H^2(\widetilde{X}_z,\mathbb{Q}))= H^2(\widetilde{X}_z,\mathbb{Q})^{\pi_1(U,z)}.\]
and if $\alpha\in H^2(\widetilde{X}_z,\mathbb{Q})^{\pi_1(U,z)}$ is a Hodge class then there is a Hodge class $\widetilde{\alpha}\in H^2(\widetilde{X},\mathbb{Q})$ such that $\widetilde{\alpha}$ restricts to $\alpha$. Since $H^2(\widetilde{X}_z,\mathbb{Q})\simeq N^1(\widetilde{X}_z)_{\mathbb{Q}}$, we see that \[\mathrm{im}(H^2(\widetilde{X},\mathbb{Q})\longrightarrow H^2(\widetilde{X}_z,\mathbb{Q})) \simeq \mathrm{im}(N^1(\widetilde{X})_{\mathbb{Q}} \longrightarrow N^1(\widetilde{X}_z)_{\mathbb{Q}})\]
for $z\in U$. In particular
\[\mathrm{im}(N^1(\widetilde{X})_{\mathbb{R}} \longrightarrow N^1(\widetilde{X}_z)_{\mathbb{R}}) \simeq
H^2(\widetilde{X}_z,\mathbb{R})^{\pi_1(U,z)}\]
for $z\in U$.\\

Now the following diagram of pull-back morphisms commutes
\[\begin{tikzcd}
N^1(X)_{\mathbb{R}} \arrow[r, "i^*"] \arrow[d, "\mu^*"]& N^1(X_z)_{\mathbb{R}} \arrow[d, "\mu_z^*"]\\
N^1(\widetilde{X})_{\mathbb{R}} \arrow[r, "\widetilde{i}^*"] & N^1(\widetilde{X}_z)_{\mathbb{R}}
\end{tikzcd}\]
Since $\mu: \widetilde{X} \longrightarrow X$ and $\mu_z: \widetilde{X}_z \longrightarrow X_z$ are resolutions of singularities for general $z\in U$, the vertical morphisms are injective. Therefore
\[\mathrm{im}(i^*)\simeq \mathrm{im}(\mu_z^*\circ i^*) = \mathrm{im}( \widetilde{i}^*\circ \mu^*)\]

Since $X$ is $\mathbb{Q}$-factorial, we have $N^1(\widetilde{X})_{\mathbb{R}} \simeq \mu^*N^1(X)_{\mathbb{R}} \oplus_i \mathbb{R}E_i$ where $E_i$ are the $\mu$-exceptional divisors. For $z\in U$ general, the restriction of a $\mu$-exceptional divisor $E_i$ to $\widetilde{X}_z$ is $\mu_z$-exceptional. In $N^1(\widetilde{X}_z)_{\mathbb{R}}$, we have $\mathrm{im}(\mu_z^*)\cap \oplus_j\mathbb{R}E_j^z= 0$ where $E_j^z$ are $\mu_z$-exceptional. Therefore
\[\mathrm{im}(\widetilde{i}^* \circ \mu^*)=  \mathrm{im}(\widetilde{i}^*)\cap \mathrm{im}(\mu_z^*).\]
Recall that we have the isomorphisms given by first Chern class $N^1(\widetilde{X}_z)_{\mathbb{R}}\simeq H^2(\widetilde{X}_z,\mathbb{R})$ and $N^1(X_z)_{\mathbb{R}}\simeq H^2(X_z,\mathbb{R})$. We know that $\mathrm{im}(\widetilde{i}^*)\simeq H^2(\widetilde{X}_z,\mathbb{R})^{\pi_1(U,z)}$ and the monodromy actions on $H^2(X_z,\mathbb{R})$ and $H^2(\widetilde{X}_z,\mathbb{R})$ commute with the pullback map $\mu_z^*$. Hence 
\[\mathrm{im}(\widetilde{i}^*)\cap \mathrm{im}(\mu_z^*)\simeq H^2(X_z,\mathbb{R})^{\pi_1(U,z)}.\]
Therefore
\[\mathrm{im}(N^1(X)_{\mathbb{R}} \longrightarrow N^1(X_z)_{\mathbb{R}})=\mathrm{im}(\widetilde{i}^*)\cap \mathrm{im}(\mu_z^*)\simeq H^2(X_z,\mathbb{R})^{\pi_1(U,z)})\]
for general $z\in U$.\\
\end{proof}

\section{Constancy on very general fibers}
Let $f:X\longrightarrow T$ be a projective morphism and $L$ is a $f$-big $\mathbb{Q}$-Cartier divisor. We denote $L_{\overline{t}}:=L|_{X_{\overline{t}}}$, the restriction to the geometric fiber of $t$. 

\begin{lem}
Let $X\longrightarrow T$ smooth projective family of varieties and $s,t \in T$ be such that $s$ is a specialization of $t$.
\begin{itemize}
    \item[(a)] $\Lambda_{\mathrm{eff}}(X_{\overline{t}})$ maps into $\Lambda_{\mathrm{eff}}(X_{\overline{s}})$ under the specialization morphism $\mathrm{sp}_{\overline{t},\overline{s}}: \mathrm{NS}_{\mathbb{R}}(X_{\overline{t}})\longrightarrow \mathrm{NS}_{\mathbb{R}}(X_{\overline{s}})$.
    \item[(b)] Suppose $a(X_{\overline{t}},L_{\overline{t}})= a(X_{\overline{s}},L_{\overline{s}})$ and $\rho(X_{\overline{t}})=\rho(X_{\overline{s}})$. Then $b(X_{\overline{t}},L_{\overline{t}})\geq b(X_{\overline{s}},L_{\overline{s}})$.
\end{itemize}
\end{lem}

\begin{proof}
(a) Let $D$ be an effective divisor in $\mathrm{NS}(X_{\overline{t}})_{\mathbb{R}}$. We may pick a discrete valuation ring $R$ with a morphism  $\phi: \mathrm{Spec}R= \{s',t' \} \longrightarrow T$ where $s'$ and $t'$ map to $s$ and $t$ respectively and $t'$ is the generic point. By Remark 2.1 we have isomorphisms $\mathrm{NS}(X_{\overline{t}}) \xrightarrow{\sim} \mathrm{NS}(X_{\overline{t'}})$ and $\mathrm{NS}(X_{\overline{s}}) \xrightarrow{\sim} \mathrm{NS}(X_{\overline{s'}})$. Therefore we may assume $T$ is the spectrum of a discrete valuation ring $R$ and $t$ is the generic point $t'$. Now $D$ is defined over a finite extension $L$ of $k(t')$. We can replace $R$ by a discrete valuation ring $R_L$ with quotient field $L$.
Then the image of $D$ under $\mathrm{Pic}(X_{t'}) \xrightarrow{\sim} \mathrm{Pic}(\phi^*X) \longrightarrow \mathrm{Pic}(X_{s'})$ is effective by semi-continuity. After passing to the algebraic closure and taking quotient by algebraic equivalence we conclude that, $\mathrm{sp}_{\overline{t},\overline{s}}$ maps $D$ to an effective divisor class.\\

(b) Since $\rho(X_{\overline{t}})=\rho(X_{\overline{s}})$,  we have an isomorphism $\mathrm{NS}(X_{\overline{t}})_{\mathbb{R}} \longrightarrow \mathrm{NS(X_{\overline{s}})}_{\mathbb{R}}$. Let $a:=a(X_{\overline{s}},L_{\overline{s}})= a(X_{\overline{t}},L_{\overline{t}})$. Note that $\mathrm{sp}_{\overline{t},\overline{s}}$ maps $K_{X_{\overline{t}}}+aL_{\overline{t}}$ to $K_{X_{\overline{s}}}+aL_{\overline{s}}$. Let $F$ be a supporting hyperplane of $\Lambda_{\mathrm{eff}}(X_{\overline{s}})$ corresponding to the minimal supporting face containing $K_{X_{\overline{s}}}+aL_{\overline{s}}$. Since $\Lambda_{\mathrm{eff}}(X_{\overline{t}}) \subset \Lambda_{\mathrm{eff}}(X_{\overline{s}})$, we see that $F$ is a supporting hyperplane of $\Lambda_{\mathrm{eff}}(X_{\overline{t}})$ containing $K_{X_{\overline{t}}}+aL_{\overline{t}}$. Therefore,
\[b(X_{\overline{s}},L_{\overline{s}})= \mathrm{codim}(F\cap \Lambda_{\mathrm{eff}}(X_{\overline{s}})) \leq  \mathrm{codim}(F\cap \Lambda_{\mathrm{eff}}(X_{\overline{t}})) \leq b(X_{\overline{t}},L_{\overline{t}}). \]
\end{proof}

\begin{lem} Let $X\longrightarrow T$ a smooth projective family. Let $\eta \in T$ be the generic point. We denote $a= a(X_{\overline{\eta}},L_{\overline{\eta}})$, $n= \rho(X_{\overline{\eta}})$ and $b=b(X_{\overline{\eta}},L_{\overline{\eta}})$.
 For $m\in \mathbb{N}$, define
\[T_m:= \{ t\in T| a(X_{\overline{t}},L_{\overline{t}})\leq a-\frac{1}{m}\}\]
\[ T_0:=\{t\in T| \rho(X_{\overline{t}})>n\}\]
and
\[T_\infty:=\{t\in T|a(X_{\overline{t}},L_{\overline{t}})=a,\rho(X_{\overline{t}})=n,b(X_{\overline{t}},L_{\overline{t}})<b\}.\]

We let $Z_T:=T_m\cup T_\infty \cup T_0$. Then
\begin{itemize}
    \item[(a)] $Z_T$ is closed under specialization.
    \item[(b)] If we base change by a morphism of schemes $g:T'\longrightarrow T$, then $Z_{T'}=g^{-1}(Z_T)$.
    
\end{itemize}
\end{lem}
\begin{proof}
 (a) Let $t\in Z_T$ and $s$ a specialization of $t$ in $T$. If $t \in T_m$ for some $m\in \mathbb{N}$, then Lemma 3.1.(a) implies that $K_{X_{\overline{s}}}+a(X_{\overline{t}},L_{\overline{t}})L_{\overline{s}} \in \Lambda_{\mathrm{eff}}(X_{\overline{s}})$. Therefore, $a(X_{\overline{t}},L_{\overline{s}}) \leq a(X_{\overline{t}},L_{\overline{t}})$ and hence $s\in T_m$. If $t \in T_0$, then by Proposition 2.1.(a), $\rho(X_{\overline{s}}) \geq \rho(X_{\overline{t}})$ and $s \in T_0$. If $t \notin T_0 \cup_m T_m$, then $\rho(X_{\overline{t}})=n$ and $a(X_{\overline{t}},L_{\overline{t}})=a$. Then Lemma 3.1.(b) implies $b(X_{\overline{s}},L_{\overline{s}}) \leq b(X_{\overline{t}},L_{\overline{t}}) < b$. Therefore $s \in T_\infty$ and $Z_T$ is closed under specialization.\\
 
 (b) This follows from the fact that the Picard number and $a$,$b$-constants are invariant under algebraically closed base extension.\\
 
 {\bf Proof of Theorem 1.1.} By passing to a resolution of singularities and using generic smoothness, we may exclude a closed subset of $T$ to assume the family $f: X \longrightarrow T$ is smooth and $T$ is affine. Since our base field $k$ is algebraically closed, we may find a subfield $k' \subset k$ which is the algebraic closure of a field  finitely generated over $\mathbb{Q}$, and there exists a finitely generated $k'$-algebra $A$ such that our family $X\longrightarrow T$ and $L$ are a base change of a family $X_A\longrightarrow \mathrm{Spec}A$ and a line bundle $L_A$ on $X_A$. Now $B=\mathrm{Spec}A$ is countable and hence $Z_B$ is a countable union, $Z_B = \cup_{b\in B} \overline{\{b\}}$ of closed subsets by Lemma 3.1.(a). Now Lemma 3.1.(b) implies that $Z_T$ is a countable union of closed subsets.
 \end{proof}
 
 \section{Family of uniruled varieties}
In this section we prove Theorem 1.2. Let $f:X \longrightarrow T$ be a projective family of uniruled varieties over an algebraically closed field $k$ of characteristic $0$ and $L$ a $f$-nef and $f$-big Cartier $\mathbb{Q}$-divisor. \\

By a standard argument using the Lefschetz principle, it is enough to prove the statement for $k=\mathbb{C}$. We will henceforth assume that $k=\mathbb{C}$. \\

We can reduce to the statement for closed points only, as follows. Let us assume that there is an open $U\subset T$ such that $b(X_t,L_t)=b$ is constant for all closed points $t\in U$. Let $s \in U$ and $Z= \overline{\{s\}}\cap U$. By applying Theorem 1.1 to the family over $Z$, we may find $F= \cup_i F_i \subset Z$ a countable union of closed subvarieties such that $b(X_{\overline{t}},L_{\overline{t}})$ is constant on $Z\setminus F$. Since $\mathbb{C}$ is uncountable, there exists a closed point $t\in Z\setminus F$. Now $s \in Z\setminus F$, since $s$ is the generic point of $Z$. Therefore, $b(X_{\overline{s}},L_{\overline{s}}) = b(X_t,L_t)=b$. Since $s \in U$ was arbitrary, we conclude that $b(X_{\overline{t}},L_{\overline{t}})=b$ for all $t\in U$. Therefore it is enough to prove the statement for closed points.

\subsection{Proof of the Theorem for closed points when $k=\mathbb{C}$.} We may replace $X$ by a resolution, and by generic smoothness, we may exclude a closed subset of the base to assume that $f: X \longrightarrow T$ is a smooth family.  By Theorem 2.5, we can shrink $T$ such that   $a(X_t,L_t)=a$ for all $t\in T$ and  $\kappa(K_{X_t}+aL_t)$ is independent of $t$. We may assume that $T$ is affine. Since $L$ is $f$-big and $f$-nef, there is an effective divisor $E$ such that for sufficiently small $s$, we have $L \sim_{\mathbb{Q}}A+sE$ where $A$ is ample. After passing to a log resolution resolving $E$, we can replace $L$ by a $\mathbb{Q}$-linearly equivalent divisor to assume that $L$ has simple normal crossings support and $(X,aL)$ is klt. \\

Since the local system $\mathcal{GN}^1(X/T)$ has finite monodromy, we can base change to a finite \'{e}tale cover of $T$ to assume that $\mathcal{GN}^1(X/T)$ has trivial monodromy.\\

If $\kappa(K_{X_t}+aL_t)=0$ then we can conclude by Proposition 2.10. Let us assume that $\kappa(K_{X_t}+aL_t)=k>0$ for all $t\in T$.\\

Since $K_X+aL$ is $f$-pseudo-effective and $aL$ is $f$-big, we may run a $K_X+aL$-MMP over $T$ to obtain a relative minimal model $\phi: X \dashrightarrow X'$. By Proposition 2.8, we may replace $T$ by an open subset to assume that the base change $\phi_t:X_t\dashrightarrow X'_t$ is a $\mathbb{Q}$-factorial minimal model and $\pi: X'_t \longrightarrow Z_t$ is the canonical model for $(X_t,aL_t)$ for all $t\in T$.\\

For $z\in Z$, we denote the image of $z$ in $T$ by $t$ and let $X'_z$ denote the fiber of $\pi: X' \longrightarrow Z$ over $z$. 

\[\begin{tikzcd}
X'_z \arrow[r] \arrow[dr] &X'_t \arrow[r], \arrow[dr, "\pi_t"]
 & X' \arrow[dr, "\pi"]& \\
&\mathrm{Spec}k(z)\arrow[r]&Z_t \arrow[r] \arrow[d, "g_t"]
& Z \arrow[d, "g"] \\
&& \mathrm{Spec}k(t) \arrow[r]
& T
\end{tikzcd}\]

Let $\mu : \widetilde{X}\longrightarrow X'$ be a resolution of singularities. We may replace $T$ by an open subset to assume that $\widetilde{X}\longrightarrow T$ is smooth. Let $\widetilde{X}_z$ be the fiber of $\widetilde{\pi}: \widetilde{X}\longrightarrow Z$ over $z\in Z$. By [Ver76, Corrolaire 5.1] we can find a Zariski open $U_Z \subset Z$ such that $\widetilde{\pi}$ is smooth over $U_Z$ and   $\widetilde{\pi}^{-1}(U_Z)\longrightarrow U_Z$ and $\pi^{-1}(U_Z)\longrightarrow U_Z$ both are topologically locally trivial fibrations (in the analytic topology). Again we may replace $T$ by a Zariski open $V \subset T$ to assume that $U_Z\longrightarrow T$ is a topologically locally trivial fibration (in the analytic topology). Let $U_t \subset Z_t$ denote the fiber of $U_Z$ over $t\in T$. \\

For all $z\in U_Z$, there is a monodromy action of $\pi_1(U_t,z)$ on  $H^2(X'_z,\mathbb{Z})$ acting by an integral matrix $M_z$ on the free part. Now for any two points $z$ and $z'$ in $U_Z$, the fundamental groups $\pi_1(U_t,z)$ and $\pi_1(U_{t'},z')$ are isomorphic, since $U_Z\longrightarrow T$ is a locally trivial fibration. Also, the cohomology groups $H^2(X'_z,\mathbb{Z})$ and $H^2(X'_{z'},\mathbb{Z})$ are isomorphic, because $\pi^{-1}(U_Z)\longrightarrow U_Z$ is a locally trivial fibration. Since the mondromy actions depend continuously on $z \in U_Z$, we see that $M_z$ is constant. Therefore the monodromy invariant subspaces have constant rank, i.e. $\mathrm{rk}H^2(X'_z,\mathbb{R})^{\pi_1(U_t,z)}$ is constant for all $z \in U_Z$.\\

By [HM07] we know that a general fiber $X'_z$ is rationally connected and has terminal singularities. Since $X'_t$ is $\mathbb{Q}$-factorial, Lemma 2.11 implies that
\[\mathrm{rk}(\mathrm{im}(N^1(X'_t)_{\mathbb{R}} \longrightarrow N^1(X'_z)_{\mathbb{R}})=
\mathrm{rk}H^2(X'_z,\mathbb{R})^{\pi_1(U_t,z)}.\]
for general $z \in U_t$. Now using Proposition 2.9.(3) we have
\[b(X_t,L_t)= \mathrm{rk}H^2(X'_z,\mathbb{R})^{\pi_1(U_t,z)}\]
for general $z\in U_t$. Since $\mathrm{rk}H^2(X'_z,\mathbb{R})^{\pi_1(U_t,z)}$ is constant for $z\in U_Z$, we may conclude that $b(X_t,L_t)$ is constant for general $t\in T$.\\

\bigskip

\noindent  Princeton University, Princeton NJ 08544-1000

{\begin{verbatim} 
akashs@math.princeton.edu
\end{verbatim}}


\begin{thebibliography}{KKMSD73}

\bibitem[BCHM10]{} C. Birkar, P. Cascini, C. D. Hacon, and J. McKernan. \emph{Existence of minimal models for
varieties of log general type}. J. Amer. Math. Soc., 23(2):405–468, 2010.

\bibitem[BDPP13]{}
S\'{e}bastien Boucksom, Jean-Pierre Demailly, Mihai Paun, and Thomas Peternell. \emph{The pseudoeffective
cone of a compact K\"{a}hler manifold and varieties of negative Kodaira dimension}. J.
Algebraic Geom., 22(2):201–248, 2013.

\bibitem[BL17]{}
T. Browning, D. Loughran,  \emph{Varieties with too many rational points} Math. Z. (2017) Online publication.

\bibitem[BM90]{} V. V. Batyrev and Yu. I. Manin. Sur le nombre des points rationnels de hauteur
born\'{e} des vari\'{e}t\'{e}s alg\'{e}briques. Math. Ann., 286(1-3):27–43, 1990.

\bibitem[Del71]{}
Pierre Deligne, Th\'{e}orie de Hodge. II, Inst. Hautes \'{E}tudes Sci. Publ. Math. 40 (1971), 5–57 (French).

\bibitem[dFH11]{}T. de Fernex and Chr. D. Hacon. \emph{Deformations of canonical pairs and Fano varieties}. J. Reine
Angew. Math., 651:97–126, 2011.

\bibitem[FMT89]{} J. Franke, Y. I. Manin, and Y. Tschinkel. Rational points of bounded height
on Fano varieties. Invent. Math., 95(2):421–435, 1989.

\bibitem[Fuj87]{} T. Fujita. \emph{On polarized manifolds whose adjoint bundles are not semipositive.}
In Algebraic geometry, Sendai, 1985, volume 10 of Adv. Stud. Pure Math.,
pages 167–178. North-Holland, Amsterdam, 1987.

\bibitem[Fuj92]{} T. Fujita. On Kodaira energy and adjoint reduction of polarized manifolds.
Manuscripta Math., 76(1):59–84, 1992.

\bibitem[HM07]{}
 Christopher D. Hacon and James McKernan. \emph{On Shokurov’s rational connectedness conjecture}.
Duke Math. J., 138(1):119–136, 2007.

\bibitem[HMX13]{}
 Christopher D. Hacon, James McKernan, and Chenyang Xu. \emph{On the birational automorphisms
of varieties of general type}. Ann. of Math. (2), 177(3):1077–1111, 2013.

\bibitem[HTT15]{}
 Brendan Hassett, Sho Tanimoto, and Yuri Tschinkel. \emph{Balanced line bundles and equivariant
compactifications of homogeneous spaces}. Int. Math. Res. Not. IMRN, (15):6375–6410, 2015.

\bibitem[KM92]{}
J. Koll\'{a}r and S. Mori, \emph{Classification of three-dimensional flips}. J. Amer. Math. Soc. 5 (1992),
no. 3, 533–703.

\bibitem[LTT16]{}
 B. Lehmann, S. Tanimoto, and Y. Tschinkel. \emph{Balanced line bundles on Fano varieties}. J. Reine
Angew. Math., 2016. Online publication.

\bibitem[LT17a]{}
B. Lehmann, S. Tanimoto, \emph{On the geometry of thin exceptional sets in Manin's Conjecture},  Duke Math. J, 2017, To appear.
  
\bibitem[LT17b]{} B. Lehmann, S. Tanimoto, \emph{Geometric Manin's Conjecture and rational curves}, 2017, Preprint 	arXiv:1702.08508 [math.AG]

\bibitem[MP12]{}
Davesh Maulik and Bjorn Poonen. \emph{N\'{e}ron-Severi groups under specialization}. Duke Math. J., 161(11):2167–2206, 2012.

\bibitem[N\'{e}r52]{}
Andr\'{e} N\'{e}ron, \emph{Probl\'{e}mes arithm\'{e}tiques et g\'{e}om\'{e}triques rattach\'{e}s \'{a} la notion de rang d’une courbe alg\'{e}brique dans un corps}, Bull. Soc. Math. France 80 (1952), 101–166 (French). 

\bibitem[Pey03]{} Emmanuel Peyre. Points de hauteur born\'{e}e, topologie ad\'{e}lique et mesures de Tamagawa. J. Th\'{e}or. Nombres Bordeaux, 15(1):319–349, 2003.

\bibitem[Rud14]{}
C\'{e}cile Le Rudulier. \emph{Points alg\'{e}briques de hauteur born\'{e}e sur une surface}, 2014.
http://cecile.lerudulier.fr/Articles/surfaces.pdf.

\bibitem[Ver76]{}
 J.-L. Verdier, \emph{Stratifications de Whitney et th\'{e}or\'{e}me de Bertini–Sard}. Invent. Math. 36
(1976), 295–312.

\bibitem[Voi03]{}
Claire Voisin, \emph{Hodge theory and complex algebraic geometry. II}, Cambridge Studies in Advanced Mathematics, vol. 77, Cambridge University Press, Cambridge, 2003 (reprinted in 2007). Translated from the French by Leila Schneps. 
\end{thebibliography}
\end{document}